%% file: 1997-040.tex
\documentstyle[12pt]{amsart}
\input epsf.sty 
\newcommand{\be}{\begin{enumerate}}
\newcommand{\ee}{\end{enumerate}}

\newcommand{\co}{{\cal C{\it o}}}
\newcommand{\cl}{{\cal L}}

\newtheorem{thm}{Theorem}[section]
\newtheorem{df}[thm]{Definition}

\numberwithin{equation}{section}

\begin{document}

\title[Order complexes of noncomplemented lattices are nonevasive]
{Order complexes of noncomplemented lattices\\are nonevasive}
              \author{Dmitry N. Kozlov}
                
     \date{\today}
\thanks{The author was supported by the Swedish Science Council grant M-PD 11292-303.}

\address {Department of Mathematics, Royal Institute of Technology,
S-100 44, Stockholm, Sweden} 
\curraddr{Mathematical Science Research Institute,
1000 Centennial Drive, Berkeley, CA 94720, U.S.A. (until May 1, 1997)}
\email{kozlov@@math.kth.se} 
 
\begin{abstract} 
  We reprove and generalize in a combinatorial way the result of 
A.~Bj\"orner, \cite[Theorem 3.3]{Bj81}, that order complexes of 
noncomplemented lattices are contractible, namely by showing that these 
simplicial complexes are in fact nonevasive, in particular collapsible.
\end{abstract}

\maketitle

                 \section{Introduction}
   
  The study of topological properties of noncomplemented lattices may
be said to have been trigged off by a paper of H.H.~Crapo, \cite{Cr}.
Namely, as a trivial corollary of his celebrated complementation
formula \cite[Theorem 3]{Cr}, one can derive that the M\"obius 
function of a noncomplemented lattice is zero. On the other hand,
the reduced Euler characteristic of an order complex of a poset,
see Definition \ref{ocd}, is equal to the M\"obius function evaluated
on that poset. Hence converting Crapo's result into algebraic
topology language gives: {\it the reduced Euler characteristic of the
order complex of a noncomplemented lattice is zero}. 

   This was later strengthened by K.~Baclawski, who showed in 
\cite[Corollary 6.3]{Ba} that {\it a finite noncomplemented lattice is
$\Bbb Z$-acyclic} (here and in the rest of the note we often identify
the poset itself with its order complex). However, Baclawski's proof
made use of the Leray spectral sequence in an intricate way and was 
far from being combinatorial.

   A further improvement was made by A.~Bj\"orner, see
\cite[Theorem 3.3]{Bj81}, where he succeeded to show that {\it a
non\-complemented lattice is contractible.} Also this proof used heavy
topology machinery, such as Nerve Lemma. Later in the joint paper of 
A.~Bj\"orner and J.W.~Walker a more precise result, known as homotopy
complementation formula, \cite[Theorem 1.1]{BW}, was given, 
enlightening the subject considerably.

   The purpose of this paper is to explain the phenomena discovered in
\cite{Bj81} in a purely combinatorial way. Namely, we shall prove that
the simplicial complex associated to a noncomplemented lattice is
nonevasive, in particular collapsible.  

   The above mentioned result of A.~Bj\"orner combined with a theorem that
can be found in \cite{Co}, allowed one to conclude that order complexes 
of noncomplemented lattices can be reduced to a point by a sequence of
collapses and anticollapses. The new fact, which can be retrieved from 
the result that we prove in this note, is that actually this reduction
can be performed avoiding anticollapses.
 
                \section{The theorem}

  Throughout the rest of this note the symbol $\Delta$ will denote
an abstract simplicial complex, that is, a subset of the set $2^{[n]}$
such that all one-element subsets of $[n]=\{1,\dots,n\}$ lie in $\Delta$
and if both $X\in\Delta$ and $Y\subset X$ then $Y\in\Delta$. We will 
adapt the usual notions of the theory of abstract simplicial complexes, 
such as $\text{lk}_{\Delta}(\sigma), \text{dl}_\Delta(\sigma)$, where
$\sigma\in\Delta$. For their description see for example \cite[(9.9)]{Bj95}.

 Let $P$ be a finite partially ordered set, shortly poset. We say that
$P$ is bounded if it contains largest and smallest elements, which we 
denote $\hat 1$ and $\hat 0$. All the posets in this note are finite and
bounded. We write $\bar P=P\setminus\{\hat 0,\hat 1\}$. 

   The main objective of this paper is to study topological objects 
associated to posets, as given by the following definition.
\begin{df} \label{ocd}
  The {\bf order complex} $\Delta(P)$ of a poset $P$ is the simplicial
complex on the vertex set $P$ whose $k$-faces are the $k$-chains in $P$.  
\end{df} 

The symbol $\cl$ will always denote a finite bounded lattice. For $x\in\cl$
we write $\co_\cl(x)$ or simply $\co(x)$ for the set of complements of 
$x$, i.e. the set $\{y\in\cl\,|\,x\land y=\hat 0$ and $x\lor y=\hat 1\}$.

\begin{df} A simplicial complex $\Delta$ is called {\bf collapsible} if 
it can be reduced to a single point by a sequence of elementary collapse
steps. An {\bf elementary collapse step} is a replacement of $\Delta$
by another simplicial complex $\Delta'=\Delta\setminus\{\sigma,\tau\}$,
where $\sigma,\tau\in\Delta$ and $\sigma$ is a proper face of exactly 
one simplex, namely $\tau$.
\end{df}

 See for example \cite[Section 11]{Bj95} for a survey of main results on 
collapsibility.

\begin{df} \label{nev}
  We say that a simplicial complex $\Delta$ on a finite number of vertices
is {\bf nonevasive} if either $\Delta$ consists of only one vertex or 
there exists a vertex $x\in \Delta$ such that both dl$_\Delta(x)$ and 
lk$_\Delta (x)$ are nonevasive.
\end{df}

 There exist several equivalent definitions of nonevasive simplicial
complexes; the one above is taken from \cite[(11.1)]{Bj95}. Originally
nonevasive complexes were defined by Kahn, Saks and Sturtevant in 
\cite{KSS}, in order to model the notion of argument complexity. A good
survey on the properties of nonevasiveness can be also found there.
For example it is proved in \cite[Proposition 1]{KSS} that nonevasiveness
implies collapsibilibity and that the implication is strict. 

 %
 %

  Now we present the main theorem of this note. 

\begin{thm} \label{main}
Let $\cl$ be a finite lattice. Let $x\in \bar\cl$ and 
let $P=\cl\setminus\co(x)$. Then $\Delta(\bar P)$ is nonevasive,
in particular it is collapsible.
\end{thm}  

{\bf Note.} The result of Theorem \ref{main} can be translated 
to algorithmic language in the following way.

Let $\cl$ be a finite lattice, $x\in\bar\cl$, $P=\cl\setminus\co_\cl(x)$.
Assume $\cal A$ is a subset of $\bar P$ which is not known in advance.
One is allowed to ask questions of type: "{\it Is $y$ in $\cal A$?\,}",
where $y\in\bar P$. Then there exists a strategy which in at most
$|\bar P|-1$ questions determines whether the set $\cal A$ is a chain in 
$\bar P$. 

See \cite[Section 3]{KSS} for more detailed description of the
interplay of topological and algorithmic properties of (non-)evasiveness. 
  
\smallskip

{\bf Proof of Theorem \ref{main}.} 
We use induction on the number of elements in the poset $\bar P$. Since 
$x\in \bar P$, we know that $\Delta(\bar P)$ is not empty. Futhermore, 
if $\bar P$ consists of only one element then $\Delta(\bar P)$ is 
a simplicial complex consisting of only one point and hence is 
nonevasive by Definition \ref{nev}.
  
  In order to show that $\Delta(\bar P)$ is nonevasive we shall find 
a suitable atom (or a coatom) $y\in P$ for which we shall prove 
the following two claims.

\smallskip

{\bf Claim 1.} {\it $P\setminus \{y\}$ is obtained from some  
lattice by taking away the complements of a certain element.}

\smallskip

{\bf Claim 2.} {\it $P_{\geq y}$ (respectively $P_{\leq y}$, 
when $y$ is a coatom) is obtained from some lattice by taking away 
the complements of some element.} 

 To conclude the result by induction we make use of the identities
$\Delta(\,\overline{P_{\geq y}}\,)=\text{lk}_{\Delta(\bar P)} (y)$,
if $y$ is an atom, and 
$\Delta(\,\overline{P\setminus\{y\}})=\text{dl}_{\Delta(\bar P)}(y)$.

\smallskip

{\bf Note.} In the rest of the proof we will always consider the case
when $y$ is an atom. The case where, $y$ is a coatom, is identical once 
the lattice is turned upside down.

 We divide the main part of the proof into two cases, depending on 
the choice of~$y$.

\smallskip

{\bf Case 1.} {\it There exists $y\in P$ such that 
\begin{enumerate}
\item $y$ is an atom (or a coatom) of $P$;
\item $y\not\leq x$ (respectively $x\not\leq y$);
\item $y\notin\co_\cl(x)$.
\end{enumerate} } 

\smallskip

{\bf Proof of Claim 1.} Since $y$ is an atom $\cl'=\cl\setminus\{y\}$ 
is a lattice. Also it is clear that $x\in\cl'$. 

If $z\in\co_{\cl'}(x)$ then $x\lor z=\hat 1$ and $x\land z=\hat 0$ in 
$\cl'$. However, the identity $x\land z=y$ is impossible in $\cl$, since 
then it would imply $y\leq x$ which we assumed to be false. Hence 
$x\land z=\hat 0$ in $\cl$ and so $z\in\co_\cl(x)$.
 Vice versa, if $z\in\co_\cl(x)$ then obviously $z\in\co_{\cl'}(x)$.

So we have shown that $\co_\cl(x)=\co_{\cl'}(x)$ and hence
$P\setminus\{y\}=\cl'\setminus\co_{\cl'}(x)$. 

\smallskip

{\bf Proof of Claim 2.} The proof is illustrated on the picture below.

 $\cl'=[y,\hat 1]$ is a lattice since it is 
an interval in the lattice $\cl$. Let $z=x\lor y$. Since $y$ is an atom
and $y\not\leq x$ we know that $x\land y=\hat 0$. On the other hand
$y\notin\co_{\cl}(x)$ so $z\neq \hat 1$.

\begin{center}
\input{pic1.tex}
\end{center}

  Let us now show that $\co_{\cl'}(z)=\co_{\cl}(x)\cap \cl'$.
It is obvious that $\co_{\cl'}(z)\supseteq\co_{\cl}(x)\cap \cl'$.
Assume that there exists $t\in\co_{\cl'}(z)\setminus\co_{\cl}(x)$. 
Then $t\land z=y$ and $t\lor z=\hat 1$.
But $$\hat 1=t\lor z=t\lor(x\lor y)=(t\lor y)\lor x=t\lor x.$$

Let $a=t\land x$. Since $t\notin\co_{\cl}(x)$ and $t\lor x=\hat 1$ 
we conclude that $a\neq \hat 0$. Now clearly 
$$z=x\lor y\geq x\geq x\land t=a$$
and
$$t\geq x\land t=a,$$
hence 
$$y=t\land z\geq a.$$

But $y$ is an atom and $a\neq\hat 0$, then $a=y$,
which contradicts the assumption $y\not\leq x$. 
 Hence $\co_{\cl'}(z)=\co_{\cl}(x)\cap \cl'$ and so 
$P_{\geq y}=\cl'\setminus\co_{\cl'}(z)$.

{\bf Case 2.} {\it For every atom (coatom) $y\in P$ either $x\geq y$
(respectively $x\leq y$) or $y\in\co(x)$.}

  Let Atom$(P)$ (Coatom$(P)$) be the set of atoms (coatoms) of $P$. 
Let $A=$ Atom$(P)\cap\co(x)$, $B=$ Coatom$(P)\cap\co(x)$. Assume 
$A\neq\emptyset$. If $a\in A$, $b\in$ Coatom$(P)$ and $b\geq a$ then $b\in B$.
Also, vice versa, if $b\in B$, $a\in$ Atom$(P)$ and $b\geq a$ then $a\in A$.
Then it is an easy exercise to see that the poset $\bar{\cal L}$ splits 
into two disjoint (though not necessarily themselves connected) parts.
The first one contains $x$ and the second one, called $S$, contains $A$, 
$B$ and all the elements comparable to the elements from these two sets.
Then $S\subseteq\co(x)$ and $\cl'=\cl\setminus S$ is still a lattice.
Thus $P=\cl'\setminus\co_{\cl'}(x)$ and since $\cl'$ has fewer elements
than $\cl$ we are done by induction.

   So we can assume that $A=B=\emptyset$, i.e. $x$ is comparable to all 
atoms and coatoms of $\cl$ and $\co_\cl(x)=\emptyset$ (so $P=\cl$).

   We can assume that $\cl$ has more than one element and hence choose 
an atom (coatom) $y\in\cl$, such that $y\neq x$. We shall now prove 
Claims 1 and 2 for this element $y$.

{\bf Proof of Claims 1 and 2.} $P\setminus\{y\}$ (respectively 
$P_{\geq y}$) is obviously a lattice and $x$ belongs to $P\setminus\{y\}$
(resp. $P_{\geq y}$). Moreover $\co_{P\setminus\{y\}}(x)=\co_{P_{\geq y}}(x)=
\emptyset$ since $x$ is less than or equal to all coatoms of 
$P\setminus\{y\}$ (resp. $P_{\geq y}$).  
\qed

{\bf Acknowledgements.} I would like to thank E.-M.~Feichtner
and G.M.~Ziegler for the careful reading of this note and A.~Bj\"orner
for pointing out the possibility of generalization of my original result.

\end{document}

%% file: pic1.tex
\begin{picture}(0,0)%
\includegraphics{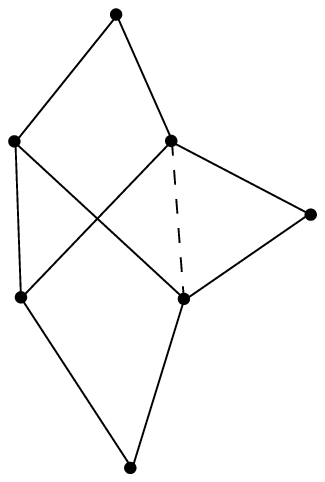}%
\end{picture}%
\setlength{\unitlength}{0.00087500in}%
\begingroup\makeatletter\ifx\SetFigFont\undefined
\def\x#1#2#3#4#5#6#7\relax{\def\x{#1#2#3#4#5#6}}%
\expandafter\x\fmtname xxxxxx\relax \def\y{splain}%
\ifx\x\y   
\gdef\SetFigFont#1#2#3{%
  \ifnum #1<17\tiny\else \ifnum #1<20\small\else
  \ifnum #1<24\normalsize\else \ifnum #1<29\large\else
  \ifnum #1<34\Large\else \ifnum #1<41\LARGE\else
     \huge\fi\fi\fi\fi\fi\fi
  \csname #3\endcsname}%
\else
\gdef\SetFigFont#1#2#3{\begingroup
  \count@#1\relax \ifnum 25<\count@\count@25\fi
  \def\x{\endgroup\@setsize\SetFigFont{#2pt}}%
  \expandafter\x
    \csname \romannumeral\the\count@ pt\expandafter\endcsname
    \csname @\romannumeral\the\count@ pt\endcsname
  \csname #3\endcsname}%
\fi
\fi\endgroup
\begin{picture}(1815,2458)(1126,-2045)
\put(1126,-1226){\makebox(0,0)[lb]{\smash{\SetFigFont{8}{9.6}{rm}$t\wedge z=y$}}}
\put(1971,329){\makebox(0,0)[lb]{\smash{\SetFigFont{8}{9.6}{rm}$\hat 1=t\vee z=t\vee x$}}}
\put(2336,-1191){\makebox(0,0)[lb]{\smash{\SetFigFont{8}{9.6}{rm}$a=x\wedge t$}}}
\put(2941,-761){\makebox(0,0)[lb]{\smash{\SetFigFont{8}{9.6}{rm}$x$}}}
\put(2311,-321){\makebox(0,0)[lb]{\smash{\SetFigFont{8}{9.6}{rm}$z=x\vee y$}}}
\put(1406,-306){\makebox(0,0)[lb]{\smash{\SetFigFont{8}{9.6}{rm}$t$}}}
\put(2021,-2021){\makebox(0,0)[lb]{\smash{\SetFigFont{8}{9.6}{rm}
$\hat 0=x\wedge y$
}}}
\end{picture}